\documentclass[a4paper,11pt]{article}

\usepackage{amsmath}
\usepackage{latexsym}

\newtheorem{theorem}{Theorem}
\newtheorem{corollary}[theorem]{Corollary}

\newtheorem{proposition}[theorem]{Proposition}
\newenvironment{proof}[1][Proof]{\textbf{#1.} }{\ \rule{0.5em}{0.5em}}

\title{ALGEBRAIC RELATIONS FOR RECURSIVE SEQUENCES}
\author{Luigi Cimmino (cimmino@na.infn.it)}
\date{}
\begin{document}
\maketitle \vspace{0.5cm}
\begin{abstract}
Through the following, we establish the conditions which allow us to
express recursive sequences of real numbers, enumerated through the
recurrence relation $a_{n+1}=A a_{n}+B a_{n-1}$, by means of
algebraic equations in two variables of degree $n\in N$. We do this,
as far as we know, like it has never been formalized before.

I'd like to precise that the work was develop without the support of
any well-know results about difference equations. Recently I
discover that the present work has been treated already, but in a
different context, at least in the two books [Car] and [Ela]. Never
the less, I found it is a good idea to publish it. It would be a
simple and ready-to-use instrument to they aren't familiar
with difference equations.
 \end{abstract}

\section{Introduction}

The famous Fibonacci's formula $a_{n+1}=a_{n}+a_{n-1}$ may be view
as the simplest case of the more general one $a_{n+1}=A a_{n}+B
a_{n-1}$. This kind of relations are usually employed in different
disciplines, like economy, pure and applied mathematics, physic,
biology and so on. So we may be interested about the conditions to be
fulfilled when one wants to express a relation like that, in the
flavour of algebraic formulas.

In what follows, we prove the statements concerning this
equivalence; then, applying our results,we show as we can compute
the domino tiling of the graph $W_{4}\times P_{n-1}$ and,
moreover, we prove rigorously that the recurrence relation
$a_{n+1}=6a_{n}-a_{n-1}$, starting with terms $a_{0}=1$ and
$a_{1}=5$, enumerates the sequence of all Pythagorean triples
$(x,y,z)$ with $y=x+1$, discussed by Mills[1].

\section{Algebraic equation of degree $n$ equals $n$-th term of the recursive sequence}

Consider the recursive sequence $\left\{a_{n}\right\}$, with $n\in
N$, obtained through the relation $a_{n+1}=A a_{n}+B a_{n-1}$
starting with $a_{0}$, $a_{1}$.

If $A$, $B$ are real number such that $A^{2}+4B>0$ and $x_{-}$,
$x_{+}$\ are the real solution of equation $x^{2}=Ax+B$, with
$x_{-}< x_{+}$, then it results $a_{n}=ax_{+}^{n}+bx_{-}^{n}$,
$\forall n\in N_{0}$, where $a=\frac{a_{1}-a_{0}x_{-}}{\sqrt{A^{2}+4B}}$ and $b=-%
\frac{a_{1}-a_{0}x_{+}}{\sqrt{A^{2}+4B}}$.\medskip

\begin{proof}
Consider the linear system of unknowns $a$ and $b$ :
\begin{eqnarray}
& &ax_{+}+bx_{-}=a_{1}\nonumber\\
& &a+b=a_{0}\nonumber
\end{eqnarray}
It admits only one solution if $x_{+}-x_{-}=\sqrt{A^{2}+4B}$
differs from zero. Now by solving the system, we obtain the
solution
\begin{displaymath}
a=\frac{a_{1}-a_{0}x_{-}}{\sqrt{A^{2}+4B}} \qquad
\textrm{and}\qquad b=- \frac{a_{1}-a_{0}x_{+}}{\sqrt{A^{2}+4B}}.
\end{displaymath}

Let $f:n\in N_{0}\rightarrow ax_{+}^{n}+bx_{-}^{n}$, so we have to
prove that $f\left( n\right) =a_{n}$, $\forall n\in N_{0}$. This
equality holds true for $n=0$, $1$. By induction, we assume
\begin{eqnarray}
& &ax_{+}^{n-1}+bx_{-}^{n-1}=a_{n-1}\nonumber\\
& &ax_{+}^{n-2}+bx_{-}^{n-2}=a_{n-2}\nonumber
\end{eqnarray}
then, by multiplying $B$ to each side of the second equation, we
have
\begin{displaymath}
aBx_{+}^{n-2}+bBx_{-}^{n-2}-Ba_{n-2}=0.
\end{displaymath}

Being $x_{\pm }^{2}=Ax_{\pm }+B$ and replacing $B=x_{\pm
}^{2}-Ax_{\pm }$,\ one writes the latter as
\begin{eqnarray}
& &\left( ax_{+}^{n-2}\right) x_{+}^{2}-\left( ax_{+}^{n-2}\right)
Ax_{+}+\left( bx_{-}^{n-2}\right) x_{-}^{2}-\left(
bx_{-}^{n-2}\right) Ax_{-}-Ba_{n-2} =\nonumber\\
& &\left( ax_{+}^{n}+bx_{-}^{n}\right) -A\left(
x_{+}^{n-1}+bx_{-}^{n-1}\right) -Ba_{n-2}=0.\nonumber
\end{eqnarray}
and substituting $ax_{+}^{n-1}+bx_{-}^{n-1}=a_{n-1}$, we have
\begin{displaymath}
ax_{+}^{n}+bx_{-}^{n}=Aa_{n-1}+Ba_{n-2}=a_{n}.
\end{displaymath}
\end{proof}\bigskip

Conversely, let us consider the following:
\begin{proposition}
Let $a$, $b$, $x_{-}$, $x_{+}$ be real number such that
$x_{-}<x_{+}$. Then, by setting $a_{n}=ax_{+}^{n}+b x_{-}^{n}$,
$\forall n\in N_{0}$, the set $S=\left\{ a_{n}\right\} _{n\in
N_{0}}$ is recursively enumerable with starting terms $a_{0}=a+b$ , \ $%
a_{1}=ax_{+}+bx_{-}$ through recurrence relation
$a_{n}=Aa_{n-1}+Ba_{n-1}$ with $A=x_{+}+x_{-}$ and $B=-x_{+}x_{-}$.
\end{proposition}

\begin{proof}
By hypotheses $a_{0}=a+b$, $a_{1}=ax_{+}+bx_{-}$ and for $n\geq
1$, one has
\begin{eqnarray}
& &a_{n}=ax_{+}^{n}+bx_{-}^{n}=\nonumber\\
& &=\left(x_{+}+x_{-}\right)
\left(ax_{+}^{n-1}+bx_{-}^{n-1}\right)
-x_{+}x_{-}\left(ax_{+}^{n-2}+bx_{-}^{n-2}\right) =\nonumber\\
& &=\left( x_{+}+x_{-}\right) a_{n-1}-x_{+}x_{-}a_{n-2}\nonumber
\end{eqnarray}
For all $n\in N_{0}$ we have $a_{n}=\left( x_{+}+x_{-}\right)
a_{n-1}-x_{+}x_{-}a_{n-2}=Aa_{n-1}+Ba_{n-2}$. \end{proof}\bigskip

We have proved the following theorem

\begin{theorem}
Let $\left\{ a_{n}\right\} _{n\in N_{0}}$ be a sequence of real
number. The following statements are equivalent
\begin{eqnarray}
& &a_{n}=Aa_{n-1}+Ba_{n-2}\\
& &a_{n}=ax_{+}^{n}+bx_{-}^{n}
\end{eqnarray} $\forall n\in N_{0}$
where $x_{+}$ and $x_{-}$ are the real solution of equation
$x^{2}=Ax+B$.
\end{theorem}

\section{A recursive sequence of natural numbers}
Before we present a first application to the graph theory, we
prove two corollaries which result useful to play with sequences
of natural numbers.

\begin{corollary}
Under the hypotheses of Proposition 1, if $%
x_{+}x_{-}<\left( x_{+}+x_{-}\right) -1$, then the sequence
$\left\{ a_{n}\right\} _{n\in N_{0}}$ is strictly increasing.
\end{corollary}

\begin{proof}
\begin{eqnarray}
& &a_{n+1} =\left( x_{+}+x_{-}\right) a_{n} -x_{+}x_{-}a_{n-1}
>\nonumber\\
& &>\left( x_{+}+x_{-}\right) a_{n} - \left[ \left(
x_{+}+x_{-}\right) -1\right] a_{n-1} =\nonumber\\
& &=\left( x_{+}+x_{-}\right) \left[ a_{n} -a_{n-1}\right]
+a_{n-1} \geq \nonumber\\
& &\geq a_{n} -a_{n-1} +a_{n-1}=a_{n}\nonumber
\end{eqnarray}
that is, $a_{n+1} >a_{n} $, $\forall n\in N_{0}$.
\end{proof}

\begin{corollary}
Under the hypotheses of Corollary 3, if $a+b$, $%
ax_{+}+bx_{-}$, $x_{+}+x_{-}$ all belong to $N$ and $x_{+}x_{-}$
$\in Z$, then $\forall n\in N_{0}$, $a_{n} \in N$.
\end{corollary}

\begin{proof}
Let $a_{0} =a+b$ and $a_{1} =ax_{+}+bx_{-}$ , both belonging to
$N$ by hypotheses. One has
\begin{eqnarray}
& &a_{2} =ax_{+}^{2}+bx_{-}^{2}=\left( x_{+}+x_{-}\right) \left(
ax_{+}+bx_{-}\right) -x_{+}x_{-}\left( a+b\right) =\nonumber\\
& &=\left( x_{+}+x_{-}\right) a_{1} -x_{+}x_{-}a_{0} \in
N\nonumber
\end{eqnarray}
So, if
\begin{eqnarray}
& &a_{n-1} =ax_{+}^{n-1}+bx_{-}^{n-1}\in N\nonumber\\
& &a_{n-2} =ax_{+}^{n-2}+bx_{-}^{n-2}\in N\nonumber
\end{eqnarray}
then
\begin{eqnarray}
& &a_{n} =ax_{+}^{n}+bx_{-}^{n}=\nonumber\\
& &=\left(x_{+}+x_{-}\right) \left(
ax_{+}^{n-1}+bx_{-}^{n-1}\right) -x_{+}x_{-}\left(
ax_{+}^{n-2}+bx_{-}^{n-2}\right) =\nonumber\\
& &=\left( x_{+}+x_{-}\right) a_{n-1} -x_{+}x_{-}a_{n-2} \in
N\nonumber
\end{eqnarray}
For all $n\in N_{0}$ we have $a_{n} \in N$.
\end{proof}\bigskip

Faase[4] proved that the number of domino tilings of the graph
$W_{4}\times P_{n-1}$, for each $n\geq 2$, can be enumerated
through recurrence relation, which is given by the product of the
$n$-th Fibonacci number and the $n$-th Pell number,{\it i.e.}
$f_{n}p_{n}= (f_{n-1}+f_{n-2})(2p_{n-1}+p_{n-2})$, as Sellers[3]
shown.

Given that the Fibonacci and Pell numbers are respectively defined
as the sequences of integers
\begin{eqnarray}
& &f_{0}=1, f_{1}=1 \textrm{ and } f_{n}=f_{n-1}+f_{n-2}\nonumber\\
& &p_{0}=1, p_{1}=2 \textrm{ and } p_{n}=2p_{n-1}+p_{n-2}\nonumber
\end{eqnarray}
applying theorem 2, we have
\begin{eqnarray}
& &f_{n} =\frac{1}{\sqrt{5}}\left[ \left(\frac{1+\sqrt{5}}{2}\right)
^{n+1}-\left(\frac{1-\sqrt{5}}{2}\right) ^{n+1}\right]\nonumber\\
& &p_{n} =\frac{1}{2\sqrt{2}}\left[ \left( 1+\sqrt{2}\right)
^{n+1}-\left( 1-\sqrt{2}\right) ^{n+1}\right]\nonumber
\end{eqnarray}
so, by multiplying them, we compute, as a function of $n$, the
number of domino tilings of these graphs, $\forall n\geq 2$.

\section{A further example}
Now, let us consider the recurrence relation
$a_{n+1}=6a_{n}-a_{n-1}$, proposed by Mills[1]. The terms of the
sequence, starting with $a_{0}=1$, $ a_{1}=5$, are the third
component of a Pythagorean triple with $y=x+1$.

By applying theorem 2, we have
\begin{displaymath}
a_{n} =\frac{1}{4}\left[ \left( \sqrt{2}+2\right) \left(
3+2\sqrt{2}\right) ^{n}-\left( \sqrt{2}-2\right) \left(
3-2\sqrt{2}\right) ^{n}\right].
\end{displaymath}
where, by mean of corollary 4, $a_{n}$ is in $N$, $\forall n\in
N_{0}$.

For all odd natural number $n$, holds true the relation $n^{2}=\left( \frac{%
n^{2}+1}{2}\right) ^{2}-\left( \frac{n^{2}-1}{2}\right) ^{2}=\left( \frac{%
n^{2}-1}{2}+1\right) ^{2}-\left( \frac{n^{2}-1}{2}\right) ^{2}$, with $%
\left( \frac{n^{2}-1}{2}\right) \in N$. Hence, the numbers
belonging to our sequence, enumerated through
$a_{k}=6a_{k-1}-a_{k-2}$ and starting with odd $a_{0}$ and
$a_{1}$, are such that $\left( \frac{a_{k}^{2}-1}{2}\right) \in
N$. Let us prove that $\forall k\in N_{0}$,
$\frac{a_{k}^{2}-1}{2}=m\left( m+1\right) $, for suitable $m\in
N$.

\begin{proof}
Let be $\frac{a_{k}^{2}-1}{2}=\left(
\frac{\sqrt{2a_{k}^{2}-1}-1}{2}\right) \left(
\frac{\sqrt{2a_{k}^{2}-1}-1}{2}+1\right) $
. Thus, we have $\frac{\sqrt{%
a_{k}^{2}-1}-1}{2}=m_{k}\in N$ \ iff $\ 2a_{k}^{2}-1=d_{k}^{2}$, $d_{k}\in N$%
.

It follows that
\begin{displaymath}
a_{k}=\frac{1}{4}\left[ \left( \sqrt{2}+2\right) \left(
3+2\sqrt{2} \right) ^{k}-\left( \sqrt{2}-2\right) \left(
3-2\sqrt{2}\right) ^{k}\right]
\end{displaymath}
then, we have
\begin{eqnarray}
& &2a_{k}^{2}-1=2\left[ \frac{1}{4}\left[ \left( \sqrt{2}+2\right)
\left( 3+2\sqrt{2}\right) ^{k}-\left( \sqrt{2}-2\right) \left(
3-2\sqrt{2}\right) ^{k}\right] \right] ^{2}-1=\nonumber\\
& &=\left[ \frac{1}{2\sqrt{2}}\left[ \left( 2+\sqrt{2}\right)
\left( 3+2\sqrt{2}\right) ^{k}-\left( 2-\sqrt{2}\right) \left(
3-2\sqrt{2}\right) ^{k}\right] \right] ^{2}=\nonumber\\
& &=\left[ \frac{1}{4}\left[ \left( 2\sqrt{2}+2\right) \left(
3+2\sqrt{2 }\right) ^{k}-\left( 2\sqrt{2}-2\right)
\left(3-2\sqrt{2}\right) ^{k}\right] \right] ^{2}\nonumber
\end{eqnarray}
that is to say
\begin{displaymath}
d_{k}=\frac{1}{2}\left[ \left( \sqrt{2}+1\right) \left( 3+2\sqrt{2}%
\right) ^{k}-\left( \sqrt{2}-1\right) \left( 3-2\sqrt{2}\right)
^{k}\right].
\end{displaymath}

By corollary 4, $\forall k\in N_{0}$, $d_{k}\in N$ and it results
a sequence of odd natural numbers for $d_{k}=6d_{k-1}-d_{k-2}$. So
it follows
$\frac{\sqrt{a_{k}^{2}-1}-1}{2}=\frac{d_{k}-1}{2}=m_{k}\in N$,
that is $ a_{k}^{2}=\left( \frac{a_{k}^{2}-1}{2}+1\right)
^{2}-\left( \frac{a_{k}^{2}-1 }{2}\right) ^{2}=\left[ m_{k}\left(
m_{k}+1\right) +1\right] ^{2}-\left[ m_{k}\left( m_{k}+1\right)
\right] ^{2}=m_{k}^{2}+\left( m_{k}+1\right) ^{2}$
.\end{proof}\bigskip

If we consider the $d_{k}$'s formula, one has
\begin{displaymath}
m_{k}=\frac{d_{k}-1}{2}=\frac{1}{4}\left[ \left( \sqrt{2}+1\right)
\left( 3+2\sqrt{2}\right) ^{k}-\left( \sqrt{2}-1\right) \left( 3-2\sqrt{2}%
\right) ^{k}-2\right]
\end{displaymath}
hence, letting $k$ vary in $N_{0}$, we have that every such
triples are enumerated by mean of the form $\left(
m_{k},m_{k}+1,a_{k}\right) $, where the $m_{k}$'s recurrence
relation is given by
\begin{eqnarray}
& &m_{k}=\frac{1}{2}\left( 6d_{k-1}-d_{k-2}-1\right)
=\frac{1}{2}\left( 6d_{k-1}-d_{k-2}-6+1+4\right)= \nonumber\\
& &= \frac{1}{2}\left[ 6\left( d_{k-1}-1\right) -\left(
d_{k-2}-1\right) +4\right] =  \nonumber\\
& &= 6\frac{\left( d_{k-1}-1\right) }{2}-\frac{\left(
d_{k-2}-1\right) }{2}+2=6m_{k-1}-m_{k-2}+2\nonumber
\end{eqnarray}
{\it i.e.} the empirical formula given by Hatch[5].

We have to prove that $\left( m_{k},m_{k}+1,a_{k}\right) $ are all
the triples with $y=x+1$. Then, we consider generic triple
$(m,m+1,c)$; from number theory there are $r$, $s\in N$ relatively
prime with different parity such that $ r>s $ and $c=r^{2}+s^{2}$.
Moreover if $m$ is even then $m=2rs$ and $ m+1=r^{2}-s^{2}=2rs+1$,
otherwise $m=r^{2}-s^{2}$ and $m+1=2rs=r^{2}-s^{2}+1$ . We obtain
two Pell's equations
\begin{eqnarray}
\textrm{m even}\qquad r^{2}-s^{2}-2rs=1& & \Rightarrow
\qquad \left(r-s\right) ^{2}-2s^{2}=1\nonumber\\
\textrm{m odd}\qquad r^{2}-s^{2}-2rs=-1& & \Rightarrow \qquad
\left( r+s\right) ^{2}-2r^{2}=1.\nonumber
\end{eqnarray}

If one assigns an index to $r$ and $s$ in each one of this cases,
then they become $\left( r_{i}-s_{i}\right) ^{2}-2s_{i}^{2}=1$ and
$\left( r_{j}+s_{j}\right) ^{2}-2r_{j}^{2}=1$. This kind of
equation admits an infinity of solution; moreover, a suitable
couple $(r_{o},s_{o})$ is either solution of the first or of the
second equation, never solution of both one. Now, we choose for
even $m_{k}$, \ let $i$ vary in the set of even natural numbers
and for odd $m_{k}$, let $j$ vary in the sets of all odd natural
numbers. We take care they are well-ordered in both one of
subsequences, therefore enumerating $m_{k}$ so that it alternates
even and odd index, we have to solve the equation
$r_{k}^{2}-s_{k}^{2}-2r_{k}s_{k}=\left( -1\right) ^{k}$, which is
$\left( r_{k}-s_{k}\right) ^{2}-2s_{k}^{2}=(-1)^{k}$. Setting
$p_{k}=r_{k}-s_{k}$ and $q_{k}=s_{k}$, we have the equation $%
p_{k}^{2}-2q_{k}^{2}=\left( -1\right) ^{k}$, whose trivial solution is $%
p_{0}=1$, $q_{0}=0$ and fundamental one $p_{1}=1$, $q_{1}=1$; from
fundamental solution we give each other solutions, since$\left( p_{1}+\sqrt{2%
}q_{1}\right) ^{k}\left( p_{1}-\sqrt{2}q_{1}\right) ^{k}=\left(
-1\right) ^{k}$, resolving system of two equations
\begin{eqnarray}
& &p_{k}+\sqrt{2}q_{k}=(1+\sqrt{2})^{k} \nonumber\\
& &p_{k}-\sqrt{2}q_{k}=(1-\sqrt{2})^{k} \nonumber
\end{eqnarray}
where clearly it holds $p_{k}+\sqrt{2}q_{k}>p_{k}-\sqrt{2}q_{k}$.

Now, subtracting and adding latter equations side to side, one
has\ respectively
\begin{eqnarray}
& &q_{k}=s_{k}=\frac{1}{4}\left[ \sqrt{2}\left( 1+\sqrt{2}\right)
^{k}- \sqrt{2}\left( 1-\sqrt{2}\right) ^{k}\right] \nonumber\\
& &p_{k}=\frac{1}{2}\left[ \left( 1+\sqrt{2}\right) ^{k}-\left(
1-\sqrt{ 2}\right) ^{k}\right] \qquad \Rightarrow \nonumber\\
& &\Rightarrow \qquad r_{k}=\frac{1}{4}\left[ \left(
\sqrt{2}+2\right) \left( 1+\sqrt{2}\right) ^{k}-\left(
\sqrt{2}-2\right) \left( 1-\sqrt{2} \right) ^{k}\right].\nonumber
\end{eqnarray}
so, stated at proposition 1, we have $r_{k}=2r_{k-1}+r_{k-2}$ and
\begin{eqnarray}
& &s_{k}=\frac{1}{4}\left[ \sqrt{2}\left( 1+\sqrt{2}\right) \left( 1+%
\sqrt{2}\right) ^{k-1}-\sqrt{2}\left( 1-\sqrt{2}\right) \left( 1-\sqrt{2}%
\right) ^{k-1}\right] =\nonumber\\
& &=\frac{1}{4}\left[ \left( \sqrt{2}+2\right) \left(
1+\sqrt{2}\right) ^{k-1}-\left( \sqrt{2}-2\right) \left(
1-\sqrt{2}\right) ^{k-1}\right] =r_{k-1} \nonumber
\end{eqnarray}
then it results
\begin{eqnarray}
& &r_{k}^{2}+s_{k}^{2}=\left(2r_{k-1}+s_{k-1}\right)
^{2}+r_{k-1}^{2}=6\left(r_{k-1}^{2}+s_{k-1}^{2}\right) -\left(
r_{k-2}^{2}+s_{k-2}^{2}\right)= \nonumber\\
& &=6c_{k-1}-c_{k-2}=c_{k}.\nonumber
\end{eqnarray}

The last result involves that the two sets $\left\{ a_{k}\right\} $ and $%
\left\{ c_{k}\right\} $, being such that $a_{0}=c_{0}=1$,
$a_{1}=c_{1}=5$ and having same recursive relation, are term by
term identical. All triples with first and second components which
mutually differ by one is in the form $\left(
2r_{n}s_{n},r_{n}^{2}-s_{n}^{2},a_{n}\right) $. To be precise, one
has
\begin{eqnarray}
& &2r_{n}s_{n}=\frac{1}{4}\left[ \left( \sqrt{2}+1\right) \left( 3+2%
\sqrt{2}\right) ^{k}-\left( \sqrt{2}-1\right) \left(
3-2\sqrt{2}\right) ^{n}-\left( -1\right) ^{n}2\right] \nonumber\\
& &r_{n}^{2}-s_{n}^{2}=\frac{1}{4}\left[ \left( \sqrt{2}+1\right)
\left( 3+2\sqrt{2}\right) ^{n}-\left( \sqrt{2}-1\right) \left( 3-2\sqrt{2}%
\right) ^{n}+\left( -1\right) ^{n}2\right] \nonumber\\
& &a_{n}=\frac{1}{4}\left[ \left( \sqrt{2}+2\right) \left( 3+2\sqrt{2}%
\right) ^{n}-\left( \sqrt{2}-2\right) \left( 3-2\sqrt{2}\right)
^{n}\right]. \nonumber
\end{eqnarray}

\section{Conclusion}
Before concluding, we note that
\begin{displaymath}
\lim_{n\to\infty}\frac{a_{n+1}}{a_{n}}=x_{+}
\end{displaymath}
which gives us the behavior of such sequences when $n$ tends to
infinity.

We have shown that the relation $a_{n+1}=A a_{n}+B a_{n-1}$ can be
express as algebraic relation, so that we can directly compute the
$n$-th term of such a sequence, without compute the preceding
terms of it. Moreover, the sequence of triples with $y=x+1$ gives
us an example who shows how we can work using results we obtained,
in the fashion of a powerful tool, to play with number theory.

\end{document}